\documentclass[11pt]{article}
\usepackage{}

\usepackage[paperwidth=8.5in,top=1.3in, bottom=1.2in]{geometry}
\usepackage{amsfonts}
\usepackage{amsthm}
\usepackage{amssymb}
\usepackage{amsmath,amsfonts}
\usepackage{mathrsfs}
\usepackage{cases}
\usepackage{latexsym,bm}
\usepackage{indentfirst}
\usepackage{color}
\usepackage{ifpdf}
\usepackage{graphicx}
\usepackage{psfrag}
\usepackage{mathtools}
\DeclarePairedDelimiterX\setc[2]{\{}{\}}{\,#1 \;\delimsize\vert\; #2\,}
\usepackage{dsfont}
\usepackage[
pdfauthor={CSY},
pdftitle={Strong generalized Halin graphs},
pdfstartview=XYZ,
bookmarks=true,
colorlinks=true,
linkcolor=blue,
urlcolor=blue,
citecolor=blue,
bookmarks=true,
linktocpage=true,
hyperindex=true
]{hyperref}
\theoremstyle{definition}

\newcommand{\xiaowuhao}{\fontsize{9pt}{\baselineskip}\selectfont}

\newtheorem{THM}{\textbf{Theorem}}[section]
\newtheorem{DEF}{\textbf{Definition}}
\newtheorem{LEM}{\textbf{Lemma}}[section]
\newtheorem{CON}{\textbf{Conjecture}}

\newcommand{\pf}{\textbf{Proof}.\quad}
\newcommand{\spf}{ \emph{Proof}.\quad}

\newtheorem{CLA}{\textbf{Claim}}[section]
\newcommand{\qqed}{\hfill $\blacksquare$\vspace{1mm}}

\newcommand{\CC}{\mathcal{H}}

\begin{document}

\title{Vizing's 2-factor Conjecture Involving Toughness and Maximum Degree  Conditions}
\vspace{1cm}
\author{Jinko Kanno$^{\dag}$ and Songling Shan$^{\ddag}$\\
{\xiaowuhao  $^{\dag}$ Louisiana Tech University, Ruston, LA \,71272 }\\
{\xiaowuhao $^{\ddag}$ Vanderbilt University, Nashville, TN\, 37240}}
\date{\today}
\maketitle

\begin{abstract}
Let $G$ be a simple graph, and let $\Delta(G)$
and $\chi'(G)$ denote the maximum degree and chromatic index of
$G$, respectively.
Vizing proved that $\chi'(G)=\Delta(G)$ or $\Delta(G)+1$.
We say $G$ is $\Delta$-critical if $\chi'(G)=\Delta+1$ and
$\chi'(H)<\chi'(G)$ for every proper subgraph $H$ of $G$.
In 1968, Vizing conjectured that if $G$ is a $\Delta$-critical graph,
then  $G$ has a 2-factor.
Let $G$ be an $n$-vertex $\Delta$-critical graph.
It was proved that if $\Delta(G)\ge n/2$, then $G$ has a 2-factor;
and that if $\Delta(G)\ge 2n/3+12$, then $G$  has a hamiltonian cycle,
and thus a 2-factor.  It is well known that every 2-tough graph with at least three
vertices has a 2-factor.
We investigate the existence of a 2-factor in a
$\Delta$-critical graph under ``moderate'' given toughness and  maximum degree conditions.
In particular, we show that  if $G$ is an  $n$-vertex $\Delta$-critical graph
with toughness at least 3/2 and with maximum degree at least $n/3$, then $G$ has a 2-factor.
In addition, we develop new techniques in proving the existence of 2-factors in graphs.
\end{abstract}

\textbf{Keywords}. Chromatic index; Critical graphs; Tutte's 2-factor theorem; Toughness

\vspace{2mm}

\section{Introduction}

In this paper, we consider only simple, undirected, and finite graphs. Let
$G$ be a graph. The notation $\Delta$  is fixed for the maximum degree of $G$ throughout the paper.
A {\it $k$-vertex\/} of $G$ is a vertex of degree
exactly $k$ in $G$. Denote by  $V_{\Delta}$ the set of $\Delta$-vertices in $G$, and by $\chi'(G)$ the chromatic index of   $G$.
The graph $G$ is called {\it critical\/}
if $\chi'(G)>\Delta$ and
$\chi'(H)<\chi'(G)$ for every proper subgraph $H$ of $G$.  It is clear that if
$G$ is critical then $G$ must be connected.
In 1965, Vizing~\cite{Vizing-2-classes} showed that a graph of maximum degree
$\Delta$ has chromatic index either $\Delta$ or $\Delta+1$.
If $\chi'(G)=\Delta$, then $G$ is said to be of class 1; otherwise, it is said to be
of class 2.
Holyer~\cite{Holyer} showed that it is NP-complete to determine whether an arbitrary graph is of class 1.
A critical graph $G$ is called {\it $\Delta$-critical\/} if
$\chi'(G)=\Delta+1$.  So $\Delta$-critical graphs are
class 2 graphs.
Motivated by the classification problem, Vizing studied critical class 2 graphs, or $\Delta$-critical graphs,
and made two well-known conjectures.

The first conjecture~\cite{vizing-ind} is on  the independence number $\alpha(G)$ of $G$, that is, the size of a maximum independent set
in $G$.

\begin{CON}[Vizing's Independence Number Conjecture]\label{ind}
Let $G$ be a $\Delta$-critical graph of order $n$.
Then $\alpha(G)\le n/2$.
\end{CON}

Furthermore, Vizing~\cite{vizing-2factor} conjectured that the following statement is true.

\begin{CON}[Vizing's 2-Factor Conjecture]\label{v2-factor}
Let $G$ be a $\Delta$-critical graph. Then $G$ contains a 2-factor.
\end{CON}

As each cycle $C$ satisfying $\alpha(C)\le |V(C)|/2$,
Conjecture~\ref{v2-factor} implies Conjecture~\ref{ind}.
%

%
%

For the Independence Number Conjecture,
Brinkmann et al.~\cite{ind1}, in 2000,
proved that if $G$   is a critical graph,
then $\alpha(G)<2n/3$; and
the upper bound is further improved when
the maximum degree is between 3 and 10.
Luo and Zhao~\cite{ind2}, in 2008, by improving the result of
Brinkmann et al., showed that if $G$ is an $n$-vertex
$\Delta$-critical graph, then
$\alpha(G)<(5\Delta-6)n/(8\Delta-6)<5n/8$ if $\Delta\ge6$.
In 2009,
Woodall~\cite{ind3} further improved the upper bound to
$3n/5$.  By restricting the problem to graphs with large maximum degrees, 
  in 2006,
 Luo and Zhao~\cite{vizing-independence-large-Delta} showed that
 Vizing's Independence Number Conjecture is true if
 $\Delta(G)\ge n/2$.

Compared to the progresses on the first Conjecture,
the progresses on Vizing's 2-Factor Conjecture has been slow.
In 2004,
Gr{\"u}newald and Steffen~\cite{vizing-2-factor-overful} established Vizing's 2-Factor Conjecture
 for graphs with the deficiency $\sum_{v\in V(G)}(\Delta(G)-d_G(v))$
small;
in particular, for overfull graphs\,(graphs
of an odd order and with the deficiency $\sum_{v\in V(G)}(\Delta(G)-d_G(v))<\Delta(G)$).
In 2012, Luo and Zhao~\cite{Vizing-2-factor-hamiltonian}
proved that if $G$ is an $n$-vertex $\Delta$-critical graph with $\Delta\ge \frac{6n}{7}$,
then $G$ contains a hamiltonian cycle, and thus a 2-factor with exactly one component.
Continuing the investigation on the existence of a hamiltonian cycle in
$\Delta$-critical graphs  with ``very large'' maximum degrees, Luo and Zhao~\cite{MR3543214} in
2016 showed that an $n$-vertex $\Delta$-critical graph with $\Delta\ge \frac{4n}{5}$ is hamiltonian.
The upper bound on $\Delta(G)$ assuring an  $n$-vertex $\Delta$-critical graph to be hamiltonian,
has been improved to $ \frac{2n}{3}+12$, respectively~\cite{chenetal2}.  Just finding 2-factors, Chen and Shan~\cite{JGT:JGT22135} proved the following result.

\begin{THM}[\cite{JGT:JGT22135}]\label{halfn}
Let $G$ be an $n$-vertex  $\Delta$-critical graph. Then $G$
has a 2-factor if $\Delta\ge n/2$.
\end{THM}

As a
measure of graph connectivity and ``resilience'' under removal of
vertices, graph toughness is a useful condition in finding factors in
graphs. To be precise, we recall the definition of toughness below.
The number of components of $G$ is denoted by $c(G)$. Let $t\ge 0$ be a
real number. The graph is said to be {\it $t$-tough\/} if $|S|\ge t\cdot
c(G-S)$ for each $S\subseteq V(G)$ with $c(G-S)\ge 2$. The {\it
toughness $\tau(G)$\/} is the largest real number $t$ for which $G$ is
$t$-tough, or is  $\infty$ if $G$ is complete. Enomoto et al.~\cite{MR785651}
proved the classic result below.

\begin{THM}[Enomoto et al.~\cite{MR785651}]\label{2-factor}
Every $k$-tough graph has a $k$-factor
if $k|V(G)|$ is even and $|V(G)|\ge k+1$.
\end{THM}

Combining the result in~Theorem~\ref{halfn} and the result in Theorem~\ref{2-factor} when
restricted to 2-factors,
one might wonder --- can we get something in between, i.e., is it possible to
find 2-factors in an $n$-vertex $\Delta$-critical graph $G$,
under the condition that $\Delta(G)<n/2$ but $\Delta(G)\ge c n$ for some
positive constant $c$, and $\tau(G)<2$ but $\tau(G)\ge d$  for some
positive constant $d$?  Particularly, we prove the following result.

\begin{THM}\label{Thm:main}
Let $G$ be an $n$-vertex  $\Delta$-critical graph. Then $G$
has a 2-factor if $\tau(G)\ge 3/2$ and $\Delta\ge n/3$.
\end{THM}

The remaining of the paper is organized as follows:
in Section 2, we recall some graph terminologies and  present several lemmas;
in Section 3, we recall Tutte's 2-factor Theorem and develop techniques
for showing the existence of 2-factors  upon applying  Tutte's 2-factor Theorem;
in the last section, we prove Theorem~\ref{Thm:main}.


%
%
%

\section{Notation and Lemmas}

Let $G$ be a graph.
For $x\in V(G)$
we denote by $d_G(x)$ the degree of $x$ in $G$.
For disjoint subsets of vertices $S$ and $T$ in $G$,
we denote by $E_G(S,T)$, the set of edges that has one
end vertex in $S$ and the other in $T$, and let $e_G(S,T)=|E_G(S,T)|$.
If  $S=\{s\}$ is a singleton,
we write $e_G(s, T)$ instead of $e_G(\{s\}, T)$.
If $H\subseteq G$ is a subgraph of $G$, and $T\subseteq V(G)$ with
$T\cap V(H)=\emptyset$, we write $E_G(H,T)$ and $e_G(H,T)$
for notational simplicity.
A {\it matching\/} in $G$ is a set of independent
edges.
If $M$ is a matching of $G$,
then let $V(M)$ denote the set of end vertices
of the edges in $M$.
For $X\subseteq V(G)$,
$M$ is said to {\it saturate\/} $X$ if $X\subseteq V(M)$.
If $G$ is a bipartite graph with partite sets $A$ and $B$, we
denote $G$ by $G[A,B]$ to emphasize the two partite sets.

To prove Theorem~\ref{Thm:main}, we present two lemmas below.

\begin{LEM}[Vizing's Adjacency Lemma]\label{vizing adjacency lemma}
Let $G$ be a $\Delta$-critical graph. Then for any edge
$xy\in E(G)$, $x$ is adjacent to at least $\Delta-d_G(y)+1$  $\Delta$-vertices
$z$ with $z\ne y$.
\end{LEM}

The following lemma
is a generalization of a
result in~\cite{vizing-independence-large-Delta}.

\begin{LEM}\label{TnoDelta}
Let $G$ be a $\Delta$-critical graph and $T$ be an independent set
in $G$.
Let $S=V(G)-T$,
and let $H=G-E(G[S])$ be the bipartite graph
with partite sets $S$ and $T$. For each $x\in S$,
let $\sigma_x$ be the number of non $\Delta$-degree neighbors of $x$ in $S$.
Assume that  there are $\delta_0$ $\Delta$-vertices
in $T$. Then for each edge
$xy\in E(H)$ with $x\in S$ and $y\in T$,
$d_H(y)\ge d_H(x)+1-\delta_0+\sigma_x$.
\end{LEM}
\pf
Let $xy\in E(H)$ with $x\in S$ and $y\in T$.
By Vizing's Adjacency Lemma,
$x$ is adjacent to at least $\Delta-d_G(y)+1$  $\Delta$-vertices
in $G$.
As $T$ has  $\delta_0$ $\Delta$-vertices,
we know $x$ is adjacent to
at least $\Delta-d_G(y)+1-\delta_0$ $\Delta$-vertices in
$S$.
Let $\sigma_x$ be the number of
all non $\Delta$-degree neighbors of $x$ in $S$.
Then, $d_H(x)+\Delta-d_G(y)+1-\delta_0+\sigma_x\le d_G(x)\le \Delta$.
By noting that $d_G(y)=d_H(y)$,
the inequality implies that $d_H(y)\ge d_H(x)+1-\delta_0+\sigma_x$.
%
\qed

\section{Tutte's 2-factor Theorem and Biased Barriers}

One of the main proof ingredients of
Theorem~\ref{Thm:main} is to apply Tutte's
2-factor Theorem under a new setting up that we develop in this section.

Let $S$ and $T$ be disjoint subsets of vertices of a graph $G$.
Let $D$ be a component of $G-(S\cup T)$.
Then $D$ is said to be an {\it odd component\/}
(resp.~{\it even component\/})
if $e_G(D, T)\equiv 1\pmod{2}$
(resp.~$e_G(D, T)\equiv 0\pmod{2}$).
Let $\CC(S, T)$ be the set of odd components of $G-(S\cup T)$
and let $h(S,T)=|\CC(S, T)|$. For $y\in T$, let
$\CC(y:S,T)=\{D\in \CC(S,T), e_G(y, D)>0\}$ and $h(y:S,T)=|\CC(y:S,T)|$.
Note that $e_G(y, V(G)-(S\cup T))\ge h(y:S,T)$.

Let
$\delta(S, T)=2|S|-2|T|+\sum_{y\in T} d_{G-S}(y)-h(S,T)$.
It is easy to see that $\delta(S, T)\equiv 0\pmod{2}$
for every $S$,~$T\subseteq V(G)$
with $S\cap T=\emptyset$.
We use the following criterion for the existence of a $2$-factor,
which is a restricted form of Tutte's $f$-factor Theorem.

\begin{LEM}[Tutte~\cite{tutte-2-factor}]\label{tutte's theorem}
A graph $G$ has a $2$-factor if and only if
$\delta(S, T)\ge 0$
for every $S$,~$T\subseteq V(G)$
with $S\cap T=\emptyset$.
\end{LEM}

An ordered pair $(S,T)$ consists of
 disjoint subsets  of vertices $S$ and $T$ in a graph $G$
is called a {\it barrier\/} if $\delta(S, T)\le -2$.
By Lemma~\ref{tutte's theorem},
if $G$ does not have a $2$-factor,
then $G$ has a barrier. We define a special barrier as below.

\begin{DEF}\label{def:biased-barrier}
Let $G$ be a graph without a 2-factor. A barrier $(S, T)$ of $G$ is called a biased barrier if among all the barriers of $G$,
\vspace{-2mm}
\begin{enumerate}
  \item $|S|$ is maximum; and
  \item subject to (1), $|T|$ is minimum.
\end{enumerate}
\end{DEF}

Properties of a minimum barrier\,(a barrier such that $|S\cup T|$ is minimum  among all the barriers of $G$)
has been established, for example, in
~\cite{AEFOS,local-CE}.
A biased barrier has similar nice properties
as given in the lemma below.

\begin{LEM}\label{biasedbarrier}
Let $G$ be a graph without  a $2$-factor,
and let $(S, T)$ be a biased barrier of $G$.
Then each of the following holds.
\begin{enumerate}
\item The set
$T$ is  independent  in $G$.
\item If $D$ is an even component
with respect to $(S, T)$,
then $e_G(T, D)=0$.
\item If $D$ is an odd component
with respect to $(S,T)$,
then for any  $y\in T$, $e_G(y, D) \le 1$.
\item If $D$ is an odd component
with respect to $(S,T)$,  then for any $x\in V(D)$,
 $e_G(x,T)\le 1$.
\end{enumerate}
\end{LEM}

\pf Let $U=V(G)-(S\cup T)$ and $z\in T$ be a vertex.
By the assumption that $(S,T)$ is a biased barrier, we know that $\delta(S, T-\{z\})\ge 0$.
So,
$$
\begin{array}{lll}
  0 &\le  \delta(S, T-\{z\})= 2|S|-2|T|+2+\sum\limits_{y\in T-\{z\}}d_{G-S}(y)-h(S,T-\{z\})& \\
  & = 2|S|-2|T|+2+\sum\limits_{y\in T}d_{G-S}(y)-e_G(z, T-\{z\})-e_G(z,U)-h(S,T-\{z\})& \\
   &\le  2|S|-2|T|+2+\sum\limits_{y\in T}d_{G-S}(y)-e_G(z, T-\{z\})-e_G(z,U)-h(S,T)+h(z:S,T)& \\
   &= \delta(S,T)+2-e_G(z, T-\{z\})-e_G(z,U)+h(z:S,T) &\\
   &\le -e_G(z, T-\{z\})-e_G(z,U)+h(z:S,T), \quad \mbox{since $\delta(S,T)\le -2$.}&
\end{array}
$$
This implies that
$$
e_G(z, T-\{z\})+e_G(z,U)-h(z:S,T)\le 0.
$$
Because $e_G(z,U)-h(z:S,T)\ge 0$ always holds, the above inequality particularly implies that
$$
e_G(z, T-\{z\})=0 \quad \mbox{for any $z\in T$} \quad \mbox{and} \quad \quad e_G(z,U)-h(z:S,T)=0.
$$
This proves statements (1)-(3).

To show (4), let $D$ be an odd component with respect to $(S,T)$
and let $x\in V(D)$
be any vertex.  Then by
the assumption that $|S|$ is maximum, we know that $\delta(S\cup \{x\}, T)\ge 0$.
So,
$$
\begin{array}{lll}
  0 &\le  \delta(S\cup \{x\}, T)= 2|S|-2|T|+2+\sum\limits_{y\in T}d_{G-(S\cup \{x\})}(y)-h(S\cup\{x\},T)& \\
  & = 2|S|-2|T|+2+\sum\limits_{y\in T}d_{G-S}(y)-e_G(x, T)-h(S\cup\{x\},T)& \\
   &\le  2|S|-2|T|+2+\sum\limits_{y\in T}d_{G-S}(y)-e_G(x, T)-(h(S,T)-1)& \\
   &= \delta(S,T)+2-e_G(x, T)+ 1&\\
   &\le -e_G(x, T)+ 1, \quad \mbox{since $\delta(S,T)\le -2$.}&
\end{array}
$$
Hence, $e_G(x, T)\le 1$.
\qed

Let $G$ be a graph without a 2-factor and let $(S,T)$
be a biased barrier of $G$. We call $(S,T)$  a
{\it good biased barrier\/} of $G$ if $h(S,T)$ is smallest among
all biased barriers of $G$.

\begin{LEM}\label{goodbb}
Let $G$ be a graph without  a $2$-factor,
and let $(S, T)$ be a good biased barrier of $G$.
For any $y\in T$,  if $h(y:S,T)\ge 2$,
then for any $D\in \CC(y: S,T)$, $|V(D)|\ge 3$.
\end{LEM}

 \pf Let $D\in  \CC(y: S,T)$ be an odd component of $G-(S\cup T)$.
 By (4) of Lemma~\ref{biasedbarrier}, $|V(D)|\ge 3$
 if $e_G(D,T)\ge 3$. So we assume that $e_G(D,T)=1$ and
 assume on the contrary that $|V(D)|\le 2$.
Let $x$ be the vertex in $D$ if $|V(D)|=1$, and
be a vertex in $D$ which is not adjacent to any vertex in $T$ if $|V(D)|=2$.
Let $z\in T$ be the vertex such that $e_G(D,z)=1$, and let
$T'=\left(T-\{z\}\right)\cup \{x\}$ and $U=V(G)-(S\cup T)$.
Let $D_z$ be the component  of $G-(S\cup T')$ which contains the
vertex $z$. Then since $e_G(z,D')=1$ for any $D'\in \CC(z:S,T)$ by (3) of Lemma~\ref{biasedbarrier},
we have that
\begin{numcases}{e_G(D_z,T)=}
\sum\limits_{D'\in \CC(z:S,T)-\{D\}}(e_G(D',T)-1)+e_G(x,z),\quad \mbox{if $|V(D)|=1$}; \nonumber \\
\sum\limits_{D'\in \CC(z:S,T)}(e_G(D',T)-1)+e_G(x,V(D)-\{x\}),\quad  \mbox{if $|V(D)|=2$}. \nonumber
\end{numcases}
Since $e_G(D',T)$ is odd for any $D'\in \CC(z:S,T)$, and $e_G(x,z)=e_G(x,V(D)-\{x\})=1$, we know that $D_z\in \CC(S,T')$ is an odd component of
$G-(S\cup T')$. Hence,  $h(S,T')=h(S,T)-h(y:S,T)+1$. So
$$
\begin{array}{lll}
   & \delta(S, T')= 2|S|-2|T|+\sum\limits_{y\in T'}d_{G-S}(y)-h(S,T') \\
   =& 2|S|-2|T|+\sum\limits_{y\in T}d_{G-S}(y)+e_G(x, V(D-x)\cup \{z\})-e_G(z, U)-h(S,T)+h(z:S,T)-1& \\
      =& \delta(S,T)+e_G(x,  V(D-x)\cup \{z\})-e_G(z, U)+h(z:S,T)-1&\\
   \le & \delta(S,T)\le -2, \quad \mbox{since $e_G(x,  V(D-x)\cup \{z\})=1, \mbox{and} \,\, e_G(z, U)\ge h(z:S,T)$.}&
\end{array}
$$
Thus, $(S, T')$  is a biased barrier. However, $h(S,T')=h(S,T)-h(y:S,T)+1\le h(S,T)-1$,
showing a contradiction to the assumption that $(S,T)$ is a good biased barrier.
\qed

 \section{Proof of Theorem~\ref{Thm:main}}

Let $G$ be an $n$-vertex  $\Delta$-critical graph such that
$\tau(G)\ge 3/2$ and $\Delta\ge n/3$. We show that $G$
has a 2-factor.

Since $G$ is 3/2-tough, $\Delta(G)\ge\delta(G)\ge 3$.
Assume
to the contrary
that $G$ does not have a $2$-factor.
Then by Tutte's 2-factor Theorem\,(Lemma~\ref{tutte's theorem}),
$G$ has a barrier.
Let $(S, T)$ be a good biased barrier of $G$.
Since $S$ and $T$ are already fixed, we simply denote $\CC(S,T)$
by $\CC$.
Let $U=V(G)-(S\cup T)$ and let
$\CC_k$ be the set of components $D$ of $G-(S\cup T)$
with $e_G(D, T)=k$.
Then we have $\CC=\bigcup_{k\ge 0}\CC_{2k+1}$.
For any  $y\in T$, let
\begin{eqnarray*}
\CC(y) & = &\setc[\big]{D\in \CC}{ e_G(y,D)=1}, \\
\CC_{1}(y) & = & \setc[\big]{D\in \CC_1 }{e_G(y,D)=1}.
\end{eqnarray*}
It is clear that $\CC_{1}(y)\subseteq \CC(y)$. Note also that $\CC(y)=\CC(y:S,T)$.
We use this notation $\CC(y)$ for simplicity  since $S$ and $T$ are already fixed.

\begin{CLA}\label{tutte}
$|T| > |S|+\sum_{k\ge 1}k|\CC_{2k+1}|$.
\end{CLA}
\spf
Since $(S, T)$ is a barrier,
\[
\begin{split}
\delta(S, T)&=2|S|-2|T|+\sum_{y\in T}d_{G-S}(y)-h(S,T)\\
&=2|S|-2|T|+\sum_{y\in T}d_{G-S}(y)-\sum_{k\ge 0} |\CC_{2k+1}| < 0.
\end{split}
\]
By Lemma~\ref{biasedbarrier}~(1) and~(2),
\[
\sum_{y\in T}d_{G-S}(y)=\sum_{y\in T}e_G(y, U)
=e_G(T, U)=\sum_{k\ge 0}(2k+1)|\CC_{2k+1}|.
\]
Therefore,
we have
\[
0 > 2|S|-2|T|+\sum_{k\ge 0}(2k+1)|\CC_{2k+1}|-\sum_{k\ge 0}|\CC_{2k+1}|,
\]
which yields
$|T| > |S|+\sum_{k\ge 1}k|\CC_{2k+1}|$.
\qed

We perform the following operations to $G$.
\begin{enumerate}
\item
Remove all even components,  and remove all components in $\CC_1$.
\item
Remove all edges in $G[S]$.
\item
For a component $D\in\CC_{2k+1}$ with $k\ge 1$
 introduce a set of $k$ independent vertices
$U^D=\{u_1^D, u_2^D,\dots, u_k^D\}$ and replace $D$
with $U^D$.
By Lemma~\ref{biasedbarrier}~(3),
$|N_G(D)\cap T|=e_G(T, D)=2k+1$.
Let $N_G(D)\cap T=\{v_0, v_1,\dots, v_{2k}\}$.
Add two new edges
$u_i^Dv_{2i-1}$
and $u_i^Dv_{2i}$ for each
$i$ with
$1\le i\le k$.
Moreover,
add one extra edge $u_1^Dv_0$.
\end{enumerate}

Let $H$ be the resulting graph, and let
$$
U^{\CC}=\bigcup_{k\ge 1}\left(\bigcup_{D\in\CC_{2k+1}}U^D\right),\quad X=S\cup U^{\CC}.
$$
By the construction,
the graph $H$ satisfies the following properties.

\begin{enumerate}
\item
$H$ is a bipartite graph with partite sets
$X$ and $T$,
\item $|U^{\CC}|=\sum_{k\ge 1} k|\CC_{2k+1}|$,
$|X|=|S|+|U^{\CC}|=|S|+\sum_{k\ge 1} k|\CC_{2k+1}|$,
and
\item
For each $k\ge 1$ and each $D\in\CC_{2k+1}$,
$d_H(u_1^D)=3$ and $d_H(u_i^D)=2$
for each $i$ with $2\le i\le k$.
\end{enumerate}
\par

We will show that
there is a matching in $H$ which saturates
$T$, which gives that $|X|=|S|+\sum_{k\ge1}k|\CC_{2k+1}|\ge |T|$,
giving  a contradiction to
Claim~\ref{tutte}.

For notation simplicity, for a set $\mathcal{D}\subseteq \CC$, let
$$
V(\mathcal{D})=\bigcup_{D\in \mathcal{D}}V(D).
$$

\begin{CLA}\label{sizeofS}
$|S|<|U^{\CC}|$.
\end{CLA}

\spf Assume on the contrary that $|S|\ge |U^{\CC}|=\sum_{k\ge 1}k|\CC_{2k+1}|$.
We may assume that $|U^{\CC}|\ge 1$. For otherwise,
since there is no edge between even component of $G-(S\cup T)$ and $T$,
and each component in $\CC_1$ is connected to a single vertex in $T$,
$c(G-S)\ge |T|$. This implies that $\tau(G)\le \frac{|S|}{|T|}<1$, giving a contradiction.

For each $D\in \CC_{2k+1}$ with $k\ge1$, let $W_D$
be a set of any $2k$ vertices in $D$
such that for each $x\in W_D$, $e_G(x,T)=1$.
Thus, $D-W_D$ is only connected  to a single vertex in $T$.
Let
$$
W=S\cup \left(\bigcup_{D\in \CC_{2k+1}, k\ge 1}W_D\right).
$$
Since $T$ is an independent set in $G$,
and each component in $G-(T\cup W)$ is connected to $S$
or only a single vertex in $T$,
 we have that $c(G-W)\ge |T|$.
So
\begin{eqnarray*}
  \tau(G) &\le & \frac{|W|}{|T|}\le  \frac{|S|+\sum_{k\ge 1}2k|\CC_{2k+1}|}{|S|+|U^{\CC}|+1}\\
   &\le &  \frac{\sum_{k\ge 1}k|\CC_{2k+1}|+\sum_{k\ge 1}2k|\CC_{2k+1}|}{\sum_{k\ge 1}k|\CC_{2k+1}|+\sum_{k\ge 1}k|\CC_{2k+1}|+1}<\frac{3}{2},
 \end{eqnarray*}
showing a contradiction to the assumption that $\tau(G)\ge 3/2$.
\qed

Because of $|T|>|S|+|U^{\CC}|$ and $|U^{\CC}|>|S|$, we get the following Claim.

\begin{CLA}\label{sizeofT}
$|T|\ge 2|S|+2$.
\end{CLA}

\begin{CLA}\label{noDelta}
$T$ contains no $\Delta$-vertex of $G$.
\end{CLA}

\spf Suppose on the contrary that there
exists $z\in T$ such that $d_G(z)=\Delta$.
We may assume that $|\CC(z)|\ge 2$.
Otherwise, $e_G(z,S)\ge \Delta-1$ and
so $|S|\ge \Delta -1$.
Hence by Claims~\ref{sizeofS} and \ref{sizeofT},
\begin{eqnarray*}
  n &\ge & |S|+|T| +|U|\\
   &\ge & 3|S|+2+|U^{\CC}|\quad (\mbox{$|U|\ge |U^{\CC}|$ by Lemma~\ref{biasedbarrier} (4)}) \\
   &\ge &  4|S|+3\ge 4\Delta-1\ge 4n/3-1,
\end{eqnarray*}
implying that $n\le 3$. This gives a contradiction
to the fact that $\Delta\ge 3$.

Hence, by Lemma~\ref{goodbb}, we have that
\begin{eqnarray*}
  n &\ge & |S|+|T| +|U|\ge e_G(z,S)+3|\CC(z)|+|T|\\
   &\ge & e_G(z,S)+3|\CC(z)|+2e_G(z,S)+2\quad (\mbox{$|T|\ge 2|S|+2\ge 2e_G(z,S)+2$} ) \\
   &= &  3(e_G(z,S)+|\CC(z)|)+2=3\Delta+2\ge n+2,
\end{eqnarray*}
showing a contradiction.
\qed

Let $D_1\in \CC$ be a component such that
$$
|V(D_1)|=\max\setc[\big]{|V(D)|}{D\in \CC},
$$
and $D_2\in \CC-\{D_1\}$ such that
$$
|V(D_2)|=\max\setc[\big]{|V(D)|}{D\in \CC-\{D_1\}}.
$$

\begin{CLA}\label{D12}
Let $D\in \CC-\{D_1,D_2\}$. Then $D$ contains no $\Delta$-vertex of $G$.
Furthermore, if $D_1$ contains a $\Delta$-vertex of $G$, then $|V(D)|\le |V(D_1)|-1$;
and if $D_2$ contains a $\Delta$-vertex of $G$, then $|V(D)|\le |V(D_2)|-2\le |V(D_1)|-2$,
 and
for any $x\in V(D)$, $d_G(x)\le \Delta-2$.
\end{CLA}

\spf Note that by the choice of $D_1$ and $D_2$, $|V(D)|\le |V(\CC)|/3$, recall here that
$V(\CC)$ is the union of vertex sets of components in $\CC$.
Since $|T|\ge 2|S|+2$ by Claim~\ref{sizeofT}, we have that
 $n \ge |S|+|T|+|V(\CC)| \ge 3|S|+2+|V(\CC)|$.  Consequently,
 $|V(\CC)|/3\le (n-2)/3-|S|$. Thus, for any $x\in V(D)$,
 \begin{eqnarray*}
  d_G(x) &\le  & |V(D)|-1+1+|S|\le |V(\CC)|/3+|S|\le (n-2)/3<\Delta.
 \end{eqnarray*}


Suppose that $D_1$ contains a $\Delta$-vertex of $G$,
and  there exists $D\in \CC-\{D_1,D_2\}$ such that $|V(D)|=|V(D_1)|$.
This implies that $|V(D_1)|=|V(D_2)|=|V(D)|$,  so $|V(D_1)|\le |V(\CC)|/3$.
Then by exactly the same argument above, we have that
for any $x\in V(D_1)$,
 \begin{eqnarray*}
  d_G(x) &\le  & |V(D)|-1+1+|S|\le |V(\CC)|/3+|S|\le (n-2)/3<\Delta.
 \end{eqnarray*}
Hence, $|V(D)|\le |V(D_1)|-1$.

Suppose now that $D_2$ contains a $\Delta$-vertex of $G$.
Since $|V(D_1)|\ge |V(D_2)|$, we then have that
$ |V(D_i)|+|S|\ge \Delta $ for $i=1,2$.
So for any $D\in \CC-\{D_1,D_2\}$,
\begin{eqnarray*}
n&\ge  & |S|+|T|+ |V(D_1)|+|V(D_2)|+|V(D)|\ge |S|+2|S|+2+|V(D_1)|+|V(D_2)|+|V(D)|\\
&\ge &|S|+|V(D_1)|+|S|+|V(D_2)|+|S|+|V(D)|+2.
 \end{eqnarray*}
Because of $ |V(D_i)|+|S|\ge \Delta $ for $i=1,2$, it follows that
$$
|S|+|V(D)|\le n-2\Delta-2\le  n/3-2\le \Delta-2.
$$
Consequently, $|V(D)|\le |V(D_2)|-2\le |V(D_1)|-2$,
 and
for any $x\in V(D)$, $d_G(x)\le \Delta-2$.
\qed

We introduce some further notation here. Let
$$
T_1=\setc[\big]{y\in T}{|\CC_1(y)|=1}, \quad \mbox{and}\quad  T_2=\setc[\big]{y\in T}{|\CC_1(y)|\ge 2}.
$$
For each component $D\in \CC_1$, let $y_D\in T$ be the
vertex such that $e_G(D, T)=e_G(D, y_D)=1$. Let
$$
\CC_{11}=\setc[\big]{D\in \CC_1}{y_D\in T_1}, \quad \mbox{and}\quad \CC_{12}=\setc[\big]{D\in \CC_1}{y_D\in T_2}(=\CC_1-\CC_{11}).
$$

\begin{CLA}\label{c-in-C12}
For each component $D\in \CC_{12}$, $|\CC(y_D)|\ge 2$. Consequently, $|V(D)|\ge 3$.
\end{CLA}

\spf Since $D\in \CC_{12}$, we have that
$|\CC(y_D)|\ge |\CC_1(y_D)|\ge 2$. The second part of the Claim is an application of Lemma~\ref{goodbb}.
\qed

Denote
\begin{eqnarray*}
 m_1 &=& |\CC_{11}|, \quad   m_2=|\CC_{12}|, \quad \mbox{and}\quad m_3=|\CC-\CC_{1}|,\\
S_1 &=& \{x\in S\,|\, x \,\text{has a non $\Delta$-degree neighbor in } V(G)-T\},\quad \text{and} \quad
S_0=S-S_1, \\
 p_y&=&|\CC_1(y)| \quad \mbox{for any $y\in T$}.
\end{eqnarray*}
Note that by the definition, if $m_2\ne 0$, then $m_2\ge 2$.


\begin{CLA}\label{extraNeighborinS}
Let $y\in T$ be a vertex. Then
  $$
  |N_G(y)\cap S|\ge
  \left\{
    \begin{array}{ll}
      |S_0|+m_1/3+m_2-1, & \hbox{if $\emptyset\ne \CC(y)\not\subseteq \{D_1,D_2\}$;} \\
      2, & \hbox{if $\CC(y)=\{D_1\}$ or $\{D_2\}$;}\\
      1, & \hbox{if $\CC(y)=\{D_1,D_2\}$.}
    \end{array}
  \right.
  $$
\end{CLA}
Moreover,  $N_G(y)\cap S\ne \emptyset$.

\spf
Since $G$ is 3/2-tough, $\delta(G)\ge 3$. As $d_G(y)=e_G(y,S)+e_G(y, V(\CC))$
and $e_G(y, V(\CC))=1$ if $y\in T_1$, so we get $e_G(y,S)\ge 2$
if $\CC(y)=\{D_1\}$ or $\{D_2\}$.

So assume that $|\CC(y)|\ge 2$. If $\CC(y)=\{D_1,D_2\}$, then $|N_G(y)\cap S|\ge 1$ as $\delta(G)\ge 3$.
Thus we assume that there exists $D\in \CC(y)-\{D_1,D_2\}$. Let $x_D$
be the neighbor of $y$ in $D$. By Claim~\ref{D12}, $x_D$ is not a $\Delta$-vertex of $G$.
Moreover, $y$ is adjacent to at least $\Delta-d_G(x_D)+1$ $\Delta$-vertices of $G$ by Vizing's
Adjacency Lemma.

Note that each component in $\CC-\CC_1$ contains at least three vertices by Lemma~\ref{biasedbarrier} (4).
So
$$
n\ge
\left\{
  \begin{array}{ll}
   |S|+|T|+|V(D_1)|+|V(D_2)|+|V(D)|+m_1+3(m_2+m_3-3), & \hbox{if $D_1, D_2, D$} \\
                                                       &\hbox{$\in \CC-\CC_{11}$;}\\
    |S|+|T|+|V(D_1)|+|V(D_2)|+|V(D)|+m_1-1+3(m_2+m_3-2), & \hbox{otherwise.}
  \end{array}
\right.
$$
Thus, because $|T|\ge 2|S|+2$ by Lemma~\ref{sizeofT}, and $|U^{\CC}|\ge |S|+1$ implying that $m_3\ge 1$,
we get that
\begin{eqnarray*}
  n &\ge & |S|+|T|+|V(D_1)|+|V(D_2)|+|V(D)|+m_1+3(m_2+m_3-3) \\
   &\ge &  3|S|+2+|V(D_1)|+|V(D_2)|+|V(D)|+m_1+3(m_2-2)\\
   &\ge &  \left\{
          \begin{array}{ll}
            3|S|+3|V(D)|+m_1+3m_2, & \hbox{if  $D_2$ contains a $\Delta$-vertex;} \\
            3|S|+3|V(D)|+m_1+3m_2-3, & \hbox{if  $D_1$  contains a $\Delta$-vertex;} \\
            3|S|+3|V(D)|+m_1+3m_2-4, & \hbox{if neither $D_1$ nor $D_2$ contains a $\Delta$-vertex.}
          \end{array}
        \right.
\end{eqnarray*}
The above bounds were obtained because of $|V(D_1)|\ge |V(D_2)|\ge |V(D)|$ and Claim~\ref{D12}.
 Thus since no component in $\CC-\{D_1, D_2\}$ containing  a $\Delta$-vertex
of $G$ by Claim~\ref{D12}, we have that
$$
|N_G(y)\cap S\cap V_{\Delta}|\ge
\left\{
  \begin{array}{ll}
    \Delta-d_G(x_D)-1, & \hbox{if  $D_2$ contains a $\Delta$-vertex;} \\
    \Delta-d_G(x_D), & \hbox{if  $D_1$  contains a $\Delta$-vertex but $D_2$ has no $\Delta$-vertex;} \\
    \Delta-d_G(x_D)+1, & \hbox{if neither $D_1$ nor $D_2$ contains a $\Delta$-vertex.}
  \end{array}
\right.
$$
Because $x_D$ is not a $\Delta$-vertex of $G$, by the definitions of $S_0$ and $S_1$,
we have that $N_G(x_D)\cap S=N_G(x_D)\cap S_1$. So
$d_G(x_D)\le |S_1|+|V(D)|$. Replacing $\Delta$ by $\frac{n}{3}$ in the above bounds on $|N_G(y)\cap S|$, and combining the  bounds on $n$,
we get that
$$
|N_G(y)\cap S|\ge
\left\{
  \begin{array}{ll}
    |S_0|+\frac{m_1}{3}+m_2-1, & \hbox{if  $D_2$ contains a $\Delta$-vertex;} \\
    |S_0|+\frac{m_1}{3}+m_2-1, & \hbox{if  $D_1$  contains a $\Delta$-vertex but $D_2$ has no $\Delta$-vertex;} \\
    |S_0|+\frac{m_1}{3}+m_2-\frac{1}{3}, & \hbox{if neither $D_1$ nor $D_2$ contains a $\Delta$-vertex.}
  \end{array}
\right.
$$

For the second part of the statement, if $\CC(y)=\emptyset$, then $N_G(y)=N_G(y)\cap S$.
So assume that $\CC(y)\ne \emptyset$. By the first part of the statement, it easily follows that $|N_G(y)\cap S|\ge 1$
unless $\CC(y)\not\subseteq \{D_1,D_2\}$. Let $D\in \CC(y)-\{D_1,D_2\}$, and let $x_D$
be the neighbor of $y$ in $D$. By Claim~\ref{D12}, $x_D$ is not a $\Delta$-vertex of $G$.
Moreover, $y$ is adjacent to at least $\Delta-d_G(x_D)+1$ $\Delta$-vertices of $G$ by Vizing's
Adjacency Lemma. Note that no component in $\CC-\{D_1, D_2\}$ contains  a $\Delta$-vertex
of $G$ by Claim~\ref{D12}. If  $D_2$ does not contain a $\Delta$-vertex of $G$,
then  $y$ is adjacent to at least $\Delta-d_G(x_D)\ge 1$ $\Delta$-vertices  which are contained
in $S$. If  $D_2$  contains a $\Delta$-vertex of $G$, then  by the second part of Claim~\ref{D12},
$d_G(x_D)\le \Delta-2$. So $y$ is adjacent to at least $\Delta-d_G(x_D)-1\ge 1$ $\Delta$-vertices  which are contained
in $S$.

The proof is finished.
\qed

If $\{D_1,D_2\}\cap \CC_1\ne \emptyset$, say $D_1\in \CC_1$, then there exists a unique  vertex $y\in T$
such that $e_G(y, D_1)=1$. We particulary name such a vertex $y$  if also $e_G(y, D_2)=1$.
\begin{equation*}
 \mbox{If $\{D_1,D_2\}\cap \CC_1\ne \emptyset$ and there exists $y\in T$ such that $\CC(y)=\{D_1,D_2\}$, we denote $y$ by $y_\omega$. }
\end{equation*}

\begin{CLA}\label{degree4vertexinT}
Let $y\in T$ be a vertex such that $|\CC(y)|\ge 2$ and $y\ne y_\omega$. Then $d_G(y)\ge 4$.
\end{CLA}

\spf Assume on the contrary that $d_G(y)=3$.  Let $D\in \CC(y)-\{D_1,D_2\}$,
and let $x_D$ be the neighbor of $y$ in $D$.
 Then $x_D$ is adjacent to at least $\Delta-3+1$ $\Delta$-vertices of
$G$ by Vizing's Adjacency Lemma.
Since $V(D)$ contains no $\Delta$-vertex of $G$ by Claim~\ref{D12},
and $T$ contains no $\Delta$-vertex of $G$ by Claim~\ref{noDelta},
we conclude that $|S|\ge |N_G(x_D)\cap S\cap V_{\Delta}|\ge \Delta-2$.
Since each $D\in \CC_{2k+1}$ contains
at least $2k+1$ vertices by Lemma~\ref{biasedbarrier} (4),
$|V(\CC)|\ge 2|U^{\CC}|$. Thus
\begin{eqnarray*}
  n &\ge & |S|+2|U^{\CC}|+|T|\ge |S|+2(|S|+1)+2|S|+2 \\
   &\ge & 5|S|+4\ge 5(\Delta-2)+4\ge 5\left(\frac{n}{3}-2\right)+4= \frac{5n}{3}-6,
 \end{eqnarray*}
implying that $n\le 9$.

By Claim~\ref{extraNeighborinS},
$N_G(y)\cap S\ne \emptyset$.
Since  $|\CC(y)|\ge 2$,  by Lemma~\ref{goodbb},
$|U|\ge |V(\CC(y))|\ge 6$. Since $|S|\ge 1$,
$|T|\ge 2|S|+2\ge 4$. Hence,
$n\ge |S|+|T|+|U|\ge 1+4+6\ge 11$, a contradiction.
\qed

\begin{CLA}\label{degreeinH}
Let $xy\in E(H)$ be an edge with $x\in X$ and $y\in T$. Then each of the following holds.
\begin{enumerate}
\item  If $x\in S_0$, then $d_H(y)+p_y\ge d_H(x)+1$.
  \item  If $x\in S_1$, then $d_H(y)+p_y\ge d_H(x)+2$.
  \item If $x\in U^{\CC}$ and $p_y=0$, then $d_H(y)\ge d_H(x)$.
  \item  If $x\in U^{\CC}$, $p_y\ge1$, and $y\ne y_\omega$, then $d_H(y)+p_y\ge d_H(x)+1$.
  \item If $x\in U^{\CC}$ and $y= y_\omega$, then $d_H(y)+p_y\ge d_H(x)$.
  \end{enumerate}
\end{CLA}

\spf Statements (1) and (2) follow from Lemma~\ref{TnoDelta} by taking $\sigma_x=0$
and 1, respectively.
The statements (3) and (5) are clear, since  $d_H(y)+p_y=d_G(y)\ge \delta(G)\ge 3$,
and $d_H(x)\le 3$ for any $x\in U^{\CC}$. Now we show statement (4).
By the assumption that $x\in U^{\CC}$ and $p_y\ge 1$, we have that $|\CC(y)|\ge 2$.
Then the statement follows by Claim~\ref{degree4vertexinT},
 since $d_H(y)+p_y=d_G(y)\ge 4$, while $d_H(x)\le 3$.
\qed

\begin{CLA}\label{Hmatching}
$H$ has a matching which saturates $T$.
\end{CLA}

\spf
Suppose on the contrary that $H$ has no matching saturating $T$. By Hall's Theorem,
there is a nonempty subset $B\subseteq T$
such that $|N_H(B)|<|B|$. Among all such subsets
with this property,  we choose $B$ with
smallest cardinality.
Let $A=N_H(B)$ and $H'=H[A\cup B]$.
Then we claim that in $H'$,  there is a matching which saturates
$A$.
Suppose this is not the case.
By Hall's Theorem again, there is a nonempty subset $A'\subseteq A$
such that $|N_{H'}(A')|<|A'|$.
Since $A'\subseteq A=N_H(B)\ne \emptyset$\,($T$ contains no isolated vertex of $H$),
$N_{H'}(A')\ne \emptyset$.  Let $B'=B-N_{H'}(A')$.
As $|B|>|A|\ge |N_{H'}(A')|>0$, $0<|B'|<|B|$.
On the other hand, we have $N_{H'}(B')=N_H(B')=A-A'$.
However, $|B'|=|B|-|N_{H'}(A')|>|A|-|N_{H'}(A')|> |A|-|A'|=|A-A'|=|N_H(B')|$,
showing a contradiction to the choice of $B$.

Let $M$ be a matching of $H'=H[A\cup B]$ which saturates $A$.
We consider two cases below.
%
%
%

{\noindent \it Case 1.}  $B\cap T_2\subseteq \{y_\omega\}$.

For any $y\in B$ with $y\ne y_\omega$, $p_y\le 1$.
Since $|B|>|A|$,
there exists $y_0\in B-V(M)$. Since $p_{y_0}\le 1$ if $y_0\ne y_\omega$, we have that $d_H(y_0)\ge 2$ if  $y_0\ne y_\omega$.
Otherwise, $d_H(y_0)\ge 1$.

Assume first that $y_\omega\not\in V(M)$.
So
applying Claim~\ref{degreeinH}, we have that
\begin{eqnarray*}
  e_G(A,B) &\le & \sum\limits_{xy\in M\atop x\in A, y\in B} d_H(x) \\
   &\le & \sum\limits_{xy\in M\atop x\in S_0\cup S_1\,\, \mbox{or}\,\,p_y=1} (d_H(y)+p_y-1)+\sum\limits_{xy\in M\atop x\in U^{\CC}, y\in B,p_y=0}d_H(y)\\
 & < &  \sum\limits_{xy\in M\atop x\in A, y\in B} d_H(y)+ d_H(y_0)\le  e_G(A,B),
\end{eqnarray*}
showing a contradiction.

Assume now that $y_\omega\in V(M)$. By the definition of $y_\omega$,   $1\le p_{y_\omega}\le 2$.
If $p_{y_\omega}=2$, then  for any edge $xy_\omega\in E(H)$, $x\in S$ and so $d_H(y_\omega)+p_{y_\omega}\ge d_H(x)+1$;
and if $p_{y_\omega}=1$, then for any edge $xy_\omega\in E(H)$, $d_H(y_\omega)+p_{y_\omega}\ge d_H(x)$.
So for any edge $xy_\omega\in E(H)$, $d_H(x)\le d_H(y_\omega)+1$.
Then
applying Claim~\ref{degreeinH}, we have that
$$
\begin{array}{lll}
  &e_G(A,B) \le \sum\limits_{xy\in M\atop x\in A, y\in B} d_H(x)& \\
   \le & \sum\limits_{xy\in M\atop x\in S_0\cup S_1\,\, \mbox{or}\,\,p_y= 1, y\ne y_\omega} (d_H(y)+p_y-1)+\left(\sum\limits_{xy\in M\atop x\in U^{\CC}, y\in B,p_y=0}d_H(y)\right)+d_H(y_\omega)+1&\\
  < & \left( \sum\limits_{xy\in M\atop x\in A, y\in B} d_H(y)\right)+d_H(y_0)\le  e_G(A,B),&
\end{array}
$$
showing a contradiction again.

{\noindent \it Case 2.}  $(B\cap T_2)-\{y_\omega\}\ne \emptyset$.

For any $y\in T_2$, $\CC_{1}(y)\subseteq \CC_{12}$. Since $|\CC_1(y)|\ge 2$
if $y\in T_2$, the assumption that $(B\cap T_2)-\{y_\omega\}\ne \emptyset$ implies that $m_2\ge 2$.
Furthermore, if $y_\omega\in T_2$, then $m_2\ge 4$.

Since $|B|>|A|$,
there exists $y_0\in B-V(M)$. Since $N_H(y)\cap S\ne \emptyset$ for any $y\in T$
by Claim~\ref{extraNeighborinS}, we have $d_H(y_0)\ge 1$.  We claim that if $y_\omega$ exists and $y_0\ne y_\omega$,
then $d_H(y_0)\ge 2$.
If $|\CC(y_0)|\le 1$, then $d_H(y_0)\ge d_G(y_0)-1\ge 2$.
So we assume that $|\CC(y_0)|\ge 2$.
If $y_\omega$ exists and $y_0\ne y_\omega$,
then by Claim~\ref{extraNeighborinS},
$d_H(y_0)\ge|N_H(y_0)\cap S|\ge |S_0|+m_1/3+m_2-1\ge m_1/3+m_2-1$.
Note that if $y_\omega\in T_2$ then $m_2\ge 4$, and if $y_\omega\not\in T_2$, then
by the definition of $y_\omega$, $m_1\ge 1$.
Thus  we have that $d_H(y_0)\ge 2$.

For any $y\in T_2-\{y_\omega\}$, $|N_H(y)\cap S|=|N_G(y)\cap S|\ge |S_0|+m_1/3+m_2-1$
by Claim~\ref{extraNeighborinS}. Thus, $|A\cap S|\ge |S_0|+m_1/3+m_2-1$.
Let $A_0=A\cap S$.
Then since $m_2\ge 2$, if $m_1\le 1$
\begin{equation}\label{sizeS0}
     2|A_0-S_0|\ge 2m_1/3+2m_2-2\ge m_1+m_2-1/3,
\end{equation}
and if $m_1\ge 2$, then
\begin{equation}\label{sizeS0}
   |A_0-S_0|\ge m_1/3+m_2-1\ge m_2-1/3.
\end{equation}

Assume first that $y_\omega\not\in V(M)$.
Applying Claim~\ref{degreeinH}, we have that
\begin{eqnarray*}
  e_G(A,B) &\le & \sum\limits_{xy\in M\atop x\in A, y\in B} d_H(x) \le \sum\limits_{xy\in M\atop x\in S_0} (d_H(y)+p_y-1)+\sum\limits_{xy\in M\atop x\in S_1} (d_H(y)+p_y-2)\\
   &&+\sum\limits_{xy\in M\atop x\in U^{\CC}, y\in B,p_y=0}d_H(y)+\sum\limits_{xy\in M\atop x\in U^{\CC}, y\in B,p_y\ge 1}(d_H(y)+p_y-1)\\
   & \le &  \left\{
            \begin{array}{ll}
              \left(\sum\limits_{xy\in M\atop x\in A, y\in B} d_H(y)\right)+  \sum\limits_{xy\in M\atop x\not\in S_1}p_y+\sum\limits_{xy\in M\atop x\in S_1}(p_y-2)& \hbox{} \\
              \left(\sum\limits_{xy\in M\atop x\in A, y\in B} d_H(y)\right)+  \sum\limits_{xy\in M\atop  p_y= 1}(p_y-1)+\sum\limits_{xy\in M\atop x\in S_1 \,\,\mbox{or}\,\, p_y\ge 2}(p_y-1) & \hbox{}
            \end{array}
          \right.\\
 & \le &  \left\{
            \begin{array}{ll}
              \left(\sum\limits_{xy\in M\atop x\in A, y\in B} d_H(y)\right)+m_1+m_2-2|A_0-S_0| & \hbox{} \\
              \left(\sum\limits_{xy\in M\atop x\in A, y\in B} d_H(y)\right)+m_2-|A_0-S_0| & \hbox{}
            \end{array}
          \right.\\
          &\le &\left(\sum\limits_{xy\in M\atop x\in A, y\in B} d_H(y)\right)+1/3\\
          &<&\left(\sum\limits_{xy\in M\atop x\in A, y\in B} d_H(y)\right)+d_H(y_0)\le e_G(A,B),
\end{eqnarray*}
showing a contradiction.

Assume now that $y_\omega\in V(M)$. By the definition of $y_\omega$,   $1\le p_{y_\omega}\le 2$.
If $p_{y_\omega}=2$, then  for any edge $xy_\omega\in E(H)$, $d_H(y_\omega)+p_{y_\omega}\ge d_H(x)+1$;
and if $p_{y_\omega}=1$, then for any edge $xy_\omega\in E(H)$, $d_H(y_\omega)+p_{y_\omega}\ge d_H(x)$.
So for any edge $xy_\omega\in E(H)$, $d_H(x)\le d_H(y_\omega)+1$.
Applying Claim~\ref{degreeinH}, we have that
\begin{eqnarray*}
  e_G(A,B) &\le & \sum\limits_{xy\in M\atop x\in A, y\in B} d_H(x) \le  \sum\limits_{xy\in M\atop x\in S_0, y\ne y_\omega} (d_H(y)+p_y-1)+\sum\limits_{xy\in M\atop x\in S_1,y\ne y_\omega} (d_H(y)+p_y-2)\\
   &&+\sum\limits_{xy\in M\atop x\in U^{\CC}, y\in B,p_y=0}d_H(y)+\sum\limits_{xy\in M\atop x\in U^{\CC}, y\in B,p_y\ge 1,y\ne y_\omega}(d_H(y)+p_y-1)+d_H(y_\omega)+1\\
 & \le &  \left\{
            \begin{array}{ll}
              \left(\sum\limits_{xy\in M\atop x\in A, y\in B} d_H(y)\right)+1+m_1+m_2-2|A_0-S_0| & \hbox{} \\
              \left(\sum\limits_{xy\in M\atop x\in A, y\in B} d_H(y)\right)+1+m_2-|A_0-S_0|& \hbox{}
            \end{array}
          \right.\\
 &\le & \left(\sum\limits_{xy\in M\atop x\in A, y\in B} d_H(y)\right)+1+1/3\\
 &<&\left(\sum\limits_{xy\in M\atop x\in A, y\in B} d_H(y)\right)+d_H(y_0)\le e_G(A,B),\quad(\mbox{$d_H(y_0)\ge 2$ in this case})
\end{eqnarray*}
showing a contradiction again.
\qed

Claim~\ref{Hmatching} gives a contradiction to Claim~\ref{tutte}.
The proof of Theorem~\ref{Thm:main} is now complete.
\qqed


\bibliographystyle{plain}
\bibliography{SSL-BIB}

\end{document}